\documentclass[12pt]{article}
\usepackage{amsfonts}
\usepackage{amssymb}
\usepackage{amsmath}
\usepackage{mathrsfs}
\usepackage{CJK}

\begin{document}
\title{  Integral-type operators from Hardy space to Bloch space on the upper half-plane}
\author{\large Ning Xu
\footnote{Corresponding author, e-mail: gx899200@126.com.}\\
\small Department of Mathematics and Science, Huai Hai Institute
of Technology,\\
\small Jiangsu, Lianyungang 222005, P. R. China}
\date{}
 \maketitle
 \vskip -1cm
\vskip 3mm

{\bf Abstract.} The boundedness of compactness of integral-type operators
from Hardy space  to Bloch space on the upper half-plane $\Pi_+=\{z\in\mathbb{C}:Imz>0\}$
are characterized.


\vskip 3mm

{\bf 1. Introduction}

Let $\Pi_+$ be the upper half-plane and $H(\Pi_+)$ be  the space of all holomorphic functions on the $\Pi_+$.
The Hardy spaces $H^2(\Pi_+)$ consists of all $f\in H(\Pi_+)$ such that
$$
\|f\|^2_{H^2(\Pi_+)}=\sup_{y>0}\int^\infty_{-\infty}|f(x+iy)|^2dx<\infty
$$
where $y=Imz$.

The Bloch space $B_\infty(\Pi_+)$ consists of all of $f\in H(\Pi_+)$ such that
$$
B(f)=\sup_{z\in{\Pi_+}}Imz|f'(z)|<\infty.
$$
A natural norm on the Bloch space can be introduced as follows
$$
\|f\|^2_{B_\infty(\Pi_+)}=|f(i)|+B(f).
$$
With this norm the $B_\infty(\Pi_+)$ becomes a Banach space.

For $f\in H(\Pi_+)$,fix $x_0\in\Pi_+$ the integral-type operators $J_g$ and $I_g$ are defined by
$$
J_g f(z)=\int_{z_0}^z f(\zeta)g'(\zeta)d\zeta, z\in\Pi_+\ \ \ \ \ (1)
$$
$$
I_g f(z)=\int_{z_0}^z f'(\zeta)g(\zeta)d\zeta, z\in\Pi_+\ \ \ \ \ (2)
$$
where $g\in H(\Pi_+)$.

The boundedness and compactness of $J_g$ and $I_g$ have been introduced on the
unit dick in[1,8,10,11,20].The operators $J_g$ and $I_g$, as well as their n-dimensional
generalizations,acting on various spaces of analytic functions, have been recently
studied,for example, in[2-7,9,19].

The importance of the operators $J_g$ and $I_g$ comes from that
$$
J_gf+I_gf=M_gf-f(z_0)g(z_0)
$$
where the multiplication operator $M_g$ is defined by $(M_gf)(z)=g(z)f(z)$.

Let $X$ and $Y$ be topological vector spaces whose topologies are given by translation-
invariant metrics $d_X$ and $d_Y$,respectively.It is said that a linear operator $T:
X\rightarrow Y$ is metrically bounded if there exists a positive constant $K$ such that
$$
d_Y(Tf,0)\leq d_X(f,0)
$$
for all $f\in X$. When $X$ and $Y$ are Banach spaces, the metrically boundedness coincides
with the boundedness of operators between Banach spaces. If we say that an operator is
bounded we will regard that it is metrically bounded. Operator $T:X\rightarrow Y$ is said
to be metrically compact if it takes every metric ball in $X$ into a relatively compact set
in $Y$.

While there is a vast literature on composition and weighted composition operators between
spaces of holomorphic functions on the unit disk,there are few papers on there operators on
spaces of functions holomorphic on the upper half-plane(see,12-18). In this paper,we consider
 the metrically boundedness and compactness of (1),(2) acting from
 $H^2(\Pi_+)$ to $B_\infty(\Pi_+)$ on the upper half-plane.

Throughout this paper, constants are denoted by $C$,they are positive and may differ
from one occurrence to the other. The notion $A\asymp B$ means that there is a
positive constant $C$ such that $C^{-1}B\leq A\leq CB$ .

{\bf 2. Auxiliary results}

{\bf Lemma 1.}Let $L=J_g(I_g)$ and $g\in H^\infty(\Pi_+)$. Then$L:H^2(\Pi_+)$ to $B_\infty(\Pi_+)$
 is metrically compact if and only if for any bounded sequence $\{f_n\}_{n\in\mathbb{C}}$ in $H^2(\Pi_+)$
  converging to zero on compact of $\Pi_+$, we have $\lim_{n\rightarrow\infty}\|Lf_n\|_{B_\infty(\Pi+)}=0$.

The proof is standard which can be found in [14] and omitted here.

{\bf Lemma 2$^{[14]}$.}Let $0<p<\infty$ and $0<a<b$. If $f\in H^2(\Pi_+)$, then
$$
\lim_{z\rightarrow\infty z\in\Gamma_{a,b}}f(z)=0,
$$
where $\Gamma_{a,b}=\{z\in\Pi_+:a\leq Imz\leq b\}$.

{\bf Lemma 3$^{[17]}$.} Let $f\in H^2(\Pi_+)$ then
$$
|f^{(n)}(z)|\leq C\frac{\|f\|_{H^2(\Pi_+)}}{(Imz)^{n+\frac{1}{2}}}, n=0,1,2  \ \ \ \ \ (3)
$$

{\bf 3.Main results}

Here we formulate and prove the main results of this paper.

{\bf Theorem 1.} Let $g$ is holomorphic function of $\Pi_+$ and $g\in H^\infty(\Pi_+)$. Then

(a) $J_g:H^2(\Pi_+)\rightarrow B_\infty(\Pi_+)$ is bounded if and only if
$$
M_1=\sup_{z\in\Pi_+}(Imz)^{\frac{1}{2}}|g'(z)|<\infty.\ \ \ \ \ \ \ \ \ (4)
$$

(b)  $I_g:H^2(\Pi_+)\rightarrow B_\infty(\Pi_+)$ is bounded if and only if
$$
M_2=\sup_{z\in\Pi_+}\frac{|g(z)|}{(Imz)^{\frac{1}{2}}}<\infty.\ \ \ \ \ \ \ \ \ (5)
$$

{\bf Proof.}(a)Assume that condition in(4) hold. Then
\begin {eqnarray*}
\|J_gf\|_{B_\infty(\Pi_+)}&=&|J_g f(i)|+\sup_{z\in\Pi_+}Imz|f(z)g'(z)|\\
&=&|\int^i_{z_0}f(\zeta)g'(\zeta)d\zeta|+\sup_{z\in\Pi_+}Imz|f(z)g'(z)|\\
&\leq&\max_{|z_0|\leq|\zeta|\leq 1,\zeta\in\Pi_+}|f(\zeta)||g(i)-g(z_0)|\\
&+&\sup_{z\in\Pi_+}(Imz)^{\frac{1}{2}}|g'(z)|\|f\|_{H^2(\Pi_+)}.\\
\end {eqnarray*}
From $f\in H^2(\Pi_+)$ and $g\in H^\infty(\Pi_+)$, it follows that $J_g:H^2(\Pi_+)\rightarrow B_\infty(\Pi_+)$ is bounded.

Now assume that $J_g:H^2(\Pi_+)\rightarrow B_\infty(\Pi_+)$ is bounded. Choose
$$
f_w(z)=\frac{(Imw)^{\frac{3}{2}}}{\Pi^{\frac{1}{2}}(z-\bar{w})^2}\ \ \ \ \ (6)
$$
It is clear that $f_w\in H^2(\Pi_+)$ and that
$\|f_w\|_{H^2(\Pi_+)}=1$. Choose $w=z$, we obtain
$$
C\geq\|J_gf_w\|_{B_\infty(\Pi_+)}\geq\sup_{w\in\Pi_+}Imw|f_w(w)g'(w)|
=\sup_{w\in\Pi_+}\frac{(Imz)^{\frac{1}{2}}|g'(z)|}{4\Pi^{\frac{1}{2}}}
$$
Therefore we obtain (4).

(b) Assume that condition in (5) hold. Then
\begin {eqnarray*}
\|I_gf\|_{B_\infty(\Pi_+)}&=&|I_gf(i)|+\sup_{z\in\Pi_+}Imz|f'(z)g(z)|\\
&=&|\int^i_{z_0}f'(\zeta)g(\zeta)d\zeta|+\sup_{z\in\Pi_+}Imz|f'(z)g(z)|\\
&\leq&|f(i)-f(z_0)|\max_{|z_0|\leq|\zeta|\leq 1,\zeta\in\Pi_+}|g(\zeta)|\\
&+&\sup_{z\in\Pi_+}(Imz)^{\frac{1}{2}}|g'(z)|\|f\|_{H^2(\Pi_+)}.\\
\end {eqnarray*}
From $f\in H^2(\Pi_+)$ and $g\in H^\infty(\Pi_+)$, it follows that $I_g:H^2(\Pi_+)\rightarrow B_\infty(\Pi_+)$ is bounded.

Now assume that $ I_g:H^2(\Pi_+)\rightarrow B_\infty(\Pi_+)$ is bounded. Choose(6) and let $w=z$,
we obtain
$$
C\geq\|I_gf_w\|_{B_\infty(\Pi_+)}\geq\sup_{w\in\Pi_+}Imw|f'_w(w)g(w)|
=\sup_{w\in\Pi_+}\frac{|g(w)|}{4\pi^{\frac{1}{2}}(Imw)^{\frac{1}{2}}}.
$$
Therefore we obtain (5).

{\bf Theorem 2.} Let $g$ is holomorphic function of $\Pi_+$ and $g\in H^\infty(\Pi_+)$.

(a) $J_g:H^2(\Pi_+)\rightarrow B_\infty(\Pi_+)$ is compact, then
$$
\lim_{r\rightarrow0}\sup_{y<r}(Imz)^{\frac{1}{2}}|g'(z)|=0.\ \ \ \ \ \ \ \ \ (7)
$$

(b)  $I_g:H^2(\Pi_+)\rightarrow B_\infty(\Pi_+)$ is compact, then
$$
\lim_{r\rightarrow0}\sup_{y<r}\frac{|g(z)|}{(Imz)^{\frac{1}{2}}}=0.\ \ \ \ \ \ \ \ \ (8)
$$
where it is understood that if $\{z:Imz<r\}$ is empty for some $r>0$, the supremum is equal to zero.

{\bf Proof.}(a) Suppose that $J_g:H^2(\Pi_+)\rightarrow B_\infty(\Pi_+)$ is compact and (7) does not
hold. Then there exists a positive number $\delta$ and a sequence $\{z_n\}_{n\in\mathbb{N}}$ in $\Pi_+$ such
that $Imz_n\rightarrow0$ and
$$
(Imz_n)^{\frac{1}{2}}|g'(z_n)|>\delta,\ \ n\in\mathbb{N}.
$$
Choose
$$
f_n(z)=\frac{(Imw_n)^{\frac{3}{2}}}{\Pi^{\frac{1}{2}}(z-\bar{w_n})^2}\ \ \ \ \ (9)
$$
and let $z_n=w_n$,then $f_n$ is norm bounded in $H^2(\Pi_+)$ and $f_n\rightarrow0$ uniformly on compacts of $\Pi_+$ as
$Imz_n\rightarrow0$. By Lemma 2 it follows that a subsequence of $\{J_gf_n\}$ tends to 0 in $B_\infty(\Pi_+)$.
On the other hand,
$$
\|J_gf\|_{B_\infty(\Pi_+)}\geq Imw_n|f_n(w_n)g'(w_n)|=\frac{1}{\sqrt{\pi}}(Imw_n)^{\frac{1}{2}}|g'(w_n)|
=\frac{\delta}{\sqrt{\pi}},
$$
which is a contradiction.

(b) Suppose that  $I_g:H^2(\Pi_+)\rightarrow B_\infty(\Pi_+)$ is compact and (8) does not hold.
Then there exists a positive number $\delta$ and a sequence $\{z_n\}_{n\in\mathbb{N}}$ in $\Pi_+$ such
that $Imz_n\rightarrow0$ and
$$
\frac{|g(z_n)|}{(Imz_n)^{\frac{1}{2}}}>\delta
$$
Choose (9) and let $z_n=w_n$, we obtain
$$
\|I_gf\|_{B_\infty(\Pi_+)}\geq Imw_n|f'_n(w_n)g(w_n)|
=\frac{1}{4\sqrt{\pi}}\frac{|g(w_n)|}{(Imw_n)^{\frac{1}{2}}}
=\frac{\delta}{4\sqrt{\pi}},
$$
which is a contradiction.

{\bf Theorem 3.} Let$g$ is holomorphic function of $\Pi_+$ ,$g\in H^\infty(\Pi_+)$
and $J_g:H^2(\Pi_+)\rightarrow B_\infty(\Pi_+)$ is
bounded. Suppose that $g\in B_\infty(\Pi_+)$. Then

(a) $J_g:H^2(\Pi_+)\rightarrow B_\infty(\Pi_+)$ is compact if
$$
\lim_{r\rightarrow0}\sup_{y<r}(Imz)^{\frac{1}{2}}|g'(z)|=0.\ \ \ \ \ \ \ \ \ (10)
$$

(b) $I_g:H^2(\Pi_+)\rightarrow B_\infty(\Pi_+)$ is compact if
$$
\lim_{r\rightarrow0}\sup_{y<r}\frac{|g(z)|}{(Imz)^{\frac{1}{2}}}=0.\ \ \ \ \ \ \ \ \ (11)
$$

{\bf Proof.}(a) Assume that (10) holds. Then for every $\epsilon>0$, there exists an $M_1>0$
such that
$$
\sup_{y<r}(Imz)^{\frac{1}{2}}|g'(z)|<\epsilon,\ \ \ whenever Imz<M_1.\ \ \ \ \ (12)
$$
Assume $\{f_n\}_{n\in\mathbb{N}}$ is a sequence in $H^2(\Pi_+)$ such that
$\sup_{n\in\mathbb{N}}\|f_n\|_{H^2(\Pi_+)}\leq M$ and $f_n\rightarrow0$
uniformly on compacts of $\Pi_+$ as $n\rightarrow\infty$. Thus for $z\in\Pi_+$,such that
$Imz<M_1$ and each $n\in\mathbb{N},$ we have
$$
Imz|f_n(z)g'(z)|\leq C(Imz)^{\frac{1}{2}}|g'(z)|\|f_n\|_{H^2(\Pi_+)}\leq\epsilon CM.\ \ \ (13)
$$

By(3), we have
$$
|f_n(z)|\leq C\frac{\|f_n\|_{H^2(\Pi_+)}}{(Imz)^{\frac{1}{2}}}\leq C\frac{M}{(Imz)^{\frac{1}{2}}}.
$$
Thus there is an $M_2>M_1$ such that
$$
|f_n(z)|<\epsilon,\ \ \ \ Whenever Imz>M_2.
$$
Hence for $z\in\Pi_+$ such that $Imz>M_2$,and each $n\in\mathbb{N}$, we have
$$
(Imz)|f_n(z)g'(z)|<\epsilon\|g\|_{B^\infty(\Pi_+)}.\ \ \ \ \ \ (14)
$$
If $M_1\leq Imz\leq M_2$, then by Lemma 2,there exists an $M_3>0$ such that
$$
|f_n(z)|<\epsilon,\ \ \ \ \ whenever |Rz|>M_3.
$$
Therefore, for each $n\in\mathbb{N}$, when $|Rz|>M_3$, we have
$$
Imz|f_n(z)g'(z)|<\epsilon\|g\|_{B^\infty(\Pi_+)}.\ \ \ \ \ (15)
$$
If $M_1\leq Imz\leq M_2$ and $|Rz|\leq M_3$, then there exists some $n_0\in\mathbb{N}$ such that
$|f_n(z)|<\epsilon$ for all $n\geq n_0$ and so
$$
Imz|f_n(z)g'(z)|<\epsilon\|g\|_{B^\infty(\Pi_+)}.\ \ \ \ \ (16)
$$
Finally, we also have
$$
|J_gf_n(i)|=|\int_{z_0}^if_n(\zeta)g'(\zeta)d\zeta|\leq
\max_{|z_0|\leq|\zeta|\leq 1,\zeta\in\Pi_+}|f_n(\zeta)||g(i)-g(z_0)|\ \ \ (17)
$$
as $n\rightarrow\infty$.

Combining (12)-(17), we have$ \|J_g f_n\|_{B_\infty(\Pi_+)}<\epsilon C$, for $n\geq n_0$ and some $C>0$.
Thus by Lemma 1,$J_g:H^2(\Pi_+)\rightarrow B_\infty(\Pi_+)$ is compact.

 \vspace{0.2cm}
\begin{center}{REFERENCES}
\end{center}

\vspace{-0.3cm}
\begin{enumerate}

\item  A.Aleman, A.G.Siskakis,An integral operator on H$^p$,Complex Variables Theory and Application,28(2)(1995)149-158.

\item  A.Aleman,J.A.Cima,An integral operator on H$^p$ and Hardy's inequality,\\J.Anal.Math.85(2001)157-176.

\item  D.C. Chang,S.Li,S.Stevi$\acute{c}$,On some integral operators on the unit polydisk and the unit ball,
Taiwanese J. Math.11(5)(2007)1251-1286.

\item  D.C. Chang, S.Stevi$\acute{c}$, Estimates of an integral operator on function spaces,
Taiwanese J. Math. 7(3)(2003)423-432.

\item  D.C. Chang, S.Stevi$\acute{c}$,The generalized Ces$\grave{a}$ro operator on the unit polydisk,
Taiwanese J. Math.7(2)(2003)293-308.

\item  Z.Hu, Extended Ces$\grave{a}$ro operators on the Bloch spaces in the unit ball of $\mathbb{C}^n$,
Acta. Math.Sci.Ser.B Engl. Ed.23(4)(2003)561-566.

\item  Z.Hu, Extended Ces$\grave{a}$ro operators on mixed norm spaces, Proc. Amer. Math. Soc. 131(7)(2003)2171-2179.

\item  S.Li,S.Stevi$\acute{c}$, Products of Volterra type operator and composition operator from $H^\infty$ and Bloch spaces
to the Zygmund space,J.Math.Anal.Appl.345\\(2008)40-52.

\item  S. Li, S.Stevi$\acute{c}$,Ces$\grave{a}$ro-type operators on some spaces of analytic functions on the unit ball,Appl.Math.Comput. 208(2009)378-388.

\item  S. Li, S.Stevi$\acute{c}$,Products of integral-type operators and composition operators between  Bloch-type spaces,
J.Math.Anal.Appl.349(2009)596-610.

\item  S. Li, S.Stevi$\acute{c}$,Volterra-type Operators on Zygmund spaces,J.Inequal.Appl.(1)\\(2005)81-88.

\item V.Matache,Compact composition operators on H$^p$ of a half-plane,
Proc. Amer. Math. Soc. 127(1999)1483-1491.

\item S.D.Sharma,A.K.Sharma,S.Ahmed,Composition operators between Hardy and Bloch-type spaces of the upper half-plane,
Bull.Korean Math.Soc.44(2007)475-482.

\item S.Stevi$\acute{c}$,Ajay K.Sharma,Weighted composition operators between Hardy and growth spaces on the upper half-plane,
Appl.Math.Comput,217(2011)4928-4934.

\item  S.Stevi$\acute{c}$,Composition operators from the weighted Bergman space to the nth weighted-type space on the upper half-plane,
Appl.Math.Comput,217\\(2010)3379-3384.

\item  S.Stevi$\acute{c}$,Composition operators from the Hardy space to the Zygmund-type space on the upper half-plane,
Abstr. Appl. Anal.,2009(2009).Article ID161528,8 pages.

\item S.Stevi$\acute{c}$,Composition operators from the Hardy space to the nth weighted-type space on the unit disk and the half-plane,
Appl.Math.Comput,215(2010)3950-3955.

\item S.Stevi$\acute{c}$,Composition operators from the Hardy space to the Zygmund-type space on the upper half-plane and the unit disk,
J.Comput.Anal.Appl.12(2)\\(2010)305-312.

\item  N.Xu, Extended Ces$\grave{a}$ro operators on $\mu$-Bloch spaces in $\mathbb{C}^n$,
J. Math. Research  and  Exposition, 29(5)(2009)913-922.

\item W.Yang,On an integral-type operator between Bloch-type spaces,Appl.Math.\\Comput,215(3)(2009)954-960.

\end{enumerate}

\end{document}